\documentclass[12pt,reqno]{article}
\usepackage{amsmath,amssymb,amsfonts}
\usepackage{mathrsfs}
\usepackage{multicol,multirow}
\usepackage{rotating}
\usepackage{appendix}
\usepackage[numbers]{natbib}
\usepackage{tabstackengine}
\usepackage{fullpage}
\usepackage{mathtools}
\usepackage{scalerel}
\usepackage{xparse}
\usepackage{amsthm}
\usepackage{amscd}
\usepackage{stmaryrd}
\usepackage{latexsym}
\usepackage{epsf}
\usepackage{breakurl}
\usepackage{hyperref}
\usepackage{caption}
\usepackage{nicematrix}
\usepackage{tabstackengine}
\usepackage{wasysym}
\usepackage{hyperref}
\usepackage[all]{xy}
\newcommand{\lah}[2]{\genfrac \lfloor \rfloor {0pt} {0} {#1} {#2}} 
\newcommand{\stone}[2]{\genfrac [ ] {0pt} {0} {#1} {#2}}
\newcommand{\stonet}[2]{\genfrac [ ] {0pt} {1} {#1} {#2}}
\newcommand{\sttwo}[2]{\genfrac \{ \} {0pt} {0} {#1} {#2}}
\newcommand{\sttwot}[2]{\genfrac \{ \} {0pt} {1} {#1} {#2}}
\newcommand{\laht}[2]{\genfrac \lfloor \rfloor {0pt} {1} {#1} {#2}} 

\newtheorem{thm}{Theorem}[section]
\newtheorem{cor}[thm]{Corollary}

\newtheorem{conj}[thm]{Conjecture}

\theoremstyle{definition}
\newtheorem{defn}[thm]{Definition}
\theoremstyle{remark}
\newtheorem{rem}[thm]{Remark}

\begin{document}
\title{Recurrence Relations for Some Integer Sequences Related to Ward Numbers}

\author{Aleks \v{Z}igon Tankosi\v{c} \\ 
Gimnazija Nova Gorica \\
Delpinova ulica 9\\
5000 Nova Gorica \\
Slovenia \\
\url{zigontaleks@gmail.com} \\\\
\textit{In memory of Professor Marko Petkov\v{s}ek}}

\date{August 14, 2025}

\maketitle

\begin{abstract}
In this paper, we give recurrence relations and identities for some integer sequences related to Ward numbers such as Ward-Lah numbers, varied Ward numbers and binomial Ward numbers. Most of the sequences are entered in the On-Line Encyclopedia of Integer Sequences. We give triangular recurrence relations, horizontal recurrence relations, generating functions and some recurrence relations of higher order obtained by using Sister Celine's general algorithm. 
\end{abstract}



\section{Introduction} \label{section1}

The \textit{Partition transformation} was defined by Luschny \cite{LuschnyBlog, LuschnyGitHub} in 2016 via SageMath code as a sequence to lower infinite triangular array transformation. According to \cite{mathoverflow}, the Partition transformation is computed via 

\begin{align*}
P_{n}^{k} (a_{1}, a_{2} \cdots) = \sum_{q \in p_{k}(n)}^{} (-1)^{q_{0}} \prod_{j=0}^{\ell(q)-1} \binom{q_{j}}{q_{j}+1} a_{j+1}^{q_{j}},
\end{align*} where the sum is over integer partitions $q$ of $n=q_0+q_1+\ldots$ (with $q_0\geqslant q_1\geqslant\ldots$) such that $q_0=k$. $\ell(q)$ is the length of $q$, $q_{\ell(q)-1}>0$ and $q_{\ell(q)}=0$.

The \textit{Ward numbers} were first studied by Ward in 1934 \cite{Ward}. The \textit{Ward numbers of the first kind}  $\left\uparrow
\begin{matrix}
n \\k
\end{matrix}
\right\downarrow$ \cite[A269940]{SloaneOEIS} are for $n \geq k \geq 1$ defined by the recurrence relation
\begin{align}
\label{ward1rec}
\left\uparrow
\begin{matrix}
n \\k
\end{matrix}
\right\downarrow = (n+k-1) \left( \left\uparrow
\begin{matrix}
n-1 \\k
\end{matrix}
\right\downarrow + \left\uparrow
\begin{matrix}
n-1 \\k-1
\end{matrix}
\right\downarrow  \right) 
\end{align} and via the Partition transformation

\begin{align}
\label{kPWard1def}
\left\uparrow
\begin{matrix}
n \\k
\end{matrix}
\right\downarrow = (-1)^{k} (n+k)^{\underline{n}} P_{n}^{k} \left(1, \frac{x}{x+1}\right),
\end{align} for $x \in \mathbb{N}$. $\mathbb{N}$ denotes the set of all positive integers without 0. $(n+k)^{\underline{n}}$ is the falling factorial, defined by

\begin{align*}
x^{\underline{n}} = \prod_{k=0}^{n-1} (x-k)
\end{align*} and \begin{align}
\label{fallingfactiroalexplicitdef}
x^{\underline{n}} = \frac{x!}{(x-n)!}.
\end{align}
A signed version of Ward numbers of the first kind also appears in \cite[pp. 152]{Jordan}.

The \textit{Ward numbers of the second kind} $	\left\updownarrow
\begin{matrix}
n \\k
\end{matrix}
\right\updownarrow$ \cite[A181996, A269939]{SloaneOEIS} are for $n \geq k \geq 1$ defined by the recurrence relation

\begin{align}
\label{ward2rec}
\left\updownarrow
\begin{matrix}
n \\k
\end{matrix}
\right\updownarrow = k 	\left\updownarrow
\begin{matrix}
n-1 \\k
\end{matrix}
\right\updownarrow + (n+k-1) 	\left\updownarrow
\begin{matrix}
n-1 \\k-1
\end{matrix}
\right\updownarrow
\end{align}
and via the Partition transformation
\begin{align}
\label{kPWard2def}
\left\updownarrow
\begin{matrix}
n \\k
\end{matrix}
\right\updownarrow = (-1)^{k} (n+k)^{\underline{n}} P_{n}^{k} \left(\frac{1}{x+1}\right)
\end{align} for $x \in \mathbb{N}_{0}$. $\mathbb{N}_{0}$ denotes the set of all positive integers including 0.  Ward numbers of the second kind also appear in \cite[pp. 172]{Jordan}.
Boundary conditions are
\begin{align*}
\left\uparrow
\begin{matrix}
n \\0
\end{matrix}
\right\downarrow = \left\uparrow
\begin{matrix}
0 \\ k
\end{matrix}
\right\downarrow = 	\left\updownarrow
\begin{matrix}
n \\0
\end{matrix}
\right\updownarrow = 	\left\updownarrow
\begin{matrix}
0 \\k
\end{matrix}
\right\updownarrow = 0
\end{align*}
and 	
\begin{align*}
\left\uparrow
\begin{matrix}
0 \\0
\end{matrix}
\right\downarrow = 	\left\updownarrow
\begin{matrix}
0 \\0
\end{matrix}
\right\updownarrow = 1
\end{align*}
with \begin{align*}
\left\uparrow
\begin{matrix}
	n \\k
\end{matrix}
\right\downarrow = 	\left\updownarrow
\begin{matrix}
	n \\k
\end{matrix}
\right\updownarrow = 0
\end{align*} for $k > n$.

In 2022 \cite[A357367]{SloaneOEIS}, Luschny defined a triangular array with a summation formula 

\begin{align}
\label{luschnywardlahdef}
T(n, k) = \sum_{m=0}^{k} (-1)^{m+k} \binom{n+k}{n+m} \binom{n+m-1}{m-1} \frac{(n+m)!}{m!}.
\end{align} However, no further details about these quantities are given.

In \cite{LuschnyGitHub}, the following triangular arrays are given 

\begin{align}
\label{t1defintro}
T_{1}^{*}(n,k) = (-1)^{k} (2n)! P_{n}^{k} \left(1, \frac{x}{x+1} \right)
\end{align} for $x \in \mathbb{N}$ \cite[A268438]{SloaneOEIS},

\begin{align}
\label{t2defintro}
T_{2}^{*}(n,k) =  (-1)^{k} (2n)! P_{n}^{k} \left(\frac{1}{x+1} \right)
\end{align} for $x \in \mathbb{N}_{0}$ \cite[A268437]{SloaneOEIS} and

\begin{align}
\label{t3def}
T_{3}^{*}(n,k) = (-1)^{k} (2n)! P_{n}^{k} (1, 1, \cdots)
\end{align} with 
\begin{align}
\label{boundary1}
T_{1}^{*}(n,0) = T_{1}^{*}(0,k) = T_{2}^{*}(n,0) = T_{2}^{*}(0,k) = T_{3}^{*}(n,0) = T_{3}^{*}(0,k)=0
\end{align} and 	
\begin{align}
\label{boundary2}
T_{1}^{*}(0,0) = T_{2}^{*}(0,0) = T_{3}^{*}(0,0) = 1.
\end{align} Also, \begin{align}
\label{boundary 2a}
T_{1}^{*}(n,k) = T_{2}^{*}(n,k) = T_{3}^{*}(n,k) = 0
\end{align} for $k > n$. For special values for sequences $T_{1}^{*}(n,k)$ and $T_{2}^{*}(n,k)$ ($k=n$ and $k=1$), see \cite[A268437, A268438]{SloaneOEIS}.

Note that in \cite{LuschnyGitHub}, definitions for sequences $T_{1}^{*}(n,k), T_{2}^{*}(n,k)$ and $T_{3}^{*}(n,k)$ are written without $(-1)^{k}$, but that does not produce sequences later given in tables of values. Sequence $T_{3}^{*}(n,k)$ is not entered in \cite{SloaneOEIS}. Luschny gave invented names for these numbers for referencing them easily: \textit{StirlingCycleStar}, \textit{StirlingSetStar} and \textit{LahStar}. We change these names in Sections \ref{section3} and \ref{section4} to varied Ward and Ward-Lah numbers.

Luschny also defined two triangular arrays A268439 and A268440 in \cite{SloaneOEIS}, which are also related to Ward numbers

\begin{align}
\label{t1defcircintro}
T_{1}^{\circ} (n, k) = (-1)^{k} \frac{(2n)!}{k!(n-k)!} P_{n}^{k}\left(1, \frac{x}{x+1}\right)
\end{align} for $x \in \mathbb{N}$ and

\begin{align}
\label{t2defcircintro}
T_{2}^{\circ} (n, k) = (-1)^{k} \frac{(2n)!}{k!(n-k)!} P_{n}^{k}\left(\frac{1}{x+1}\right)
\end{align}
for  $x \in \mathbb{N}_{0}$ with \begin{align} \label{boundary3}
T_{1}^{\circ}(n,0) = T_{1}^{\circ}(0,k) = T_{2}^{\circ}(n,0) = T_{2}^{\circ}(0,k) = 0
\end{align}
and 	
\begin{align}
\label{boundary4}
T_{1}^{\circ}(0,0) = T_{2}^{\circ}(0,0) = 1.
\end{align} Also, \begin{align}
\label{boundary 4a}
T_{1}^{\circ}(n,k) = T_{2}^{\circ}(n,k) = 0
\end{align} for $k > n$. For special values ($k=n$ and $k=1$), see \cite[A268439, A268440]{SloaneOEIS}. We name these numbers in Section \ref{section5} as binomial Ward numbers. 

The \textit{Lah numbers} $\laht{n}{k}$ were introduced by Lah in 1954 and for $n, k \ge 1$ they satisfy an explicit formula
\begin{align}
\label{lahexplicit}
\lah{n}{k} = \frac{n!}{k!} \binom{n-1}{k-1}
\end{align} (see, for instance, \cite{petkovsekpisanski}). 
According to \cite{LuschnyBlog}, the Lah numbers also satisfy 

\begin{align}
\label{lahpartitiondef}
\lah{n}{k} = (-1)^{k} \frac{n!}{k!} P_{n}^{k} (1, 1, \cdots).
\end{align}
From this relation, we also get 

\begin{align}
P_{n}^{k} (1, 1, \cdots) = (-1)^{k} \binom{n-1}{k-1}. 
\label{pnk1,1,1}
\end{align}

The paper is organized as follows. In Section \ref{section2}, we give two triangular recurrence relations: a horizontal recurrence relation and a recurrence relation of order 3 for triangular array, defined by (\ref{luschnywardlahdef}). We also give an exponential generating function. We call these numbers Ward-Lah numbers. In Section \ref{section3}, we give two triangular recurrence relations for two integer sequences, defined by (\ref{t1defintro}) and (\ref{t2defintro}). We call these quantities varied Ward numbers of both kinds. In Section \ref{section4}, we give a triangular, a horizontal recurrence relation and a generating function for an integer sequence, defined by (\ref{t3def}). We call these numbers varied Ward-Lah numbers. In Section \ref{section5}, we give recurrence relations for two integer sequences, defined by equations (\ref{t1defcircintro}) and (\ref{t2defcircintro}). We call these quantities binomial Ward numbers. We also give two conjectured relations between central Stirling and binomial Ward numbers. Finally, in Section \ref{section6}, we analogously define binomial Ward-Lah numbers, study their recurrence relations and give a relation connecting Ward-Lah numbers and central Lah numbers. 

\section{Ward-Lah Numbers}  \label{section2}

\subsection{Definition and Explicit Formula}

We start by introducing the \textit{Ward-Lah numbers} $\genfrac\lvert\rvert{0pt}{}{n} {k} $. These numbers were introduced by Luschny in \cite[A357367]{SloaneOEIS}. 

\begin{defn}
\label{wardLahdef}
The Ward-Lah numbers are for $n, k \in \mathbb{N}_{0}$ and $n \geq k$ defined via the Partition transformation as
\begin{align}
\label{wardlahdefformula}
\genfrac\lvert\rvert{0pt}{}{n} {k} = (-1)^{k} (n+k)^{\underline{n}} P_{n}^{k} (1, 1, \cdots)
\end{align}
with $\genfrac\lvert\rvert{0pt}{}{0} {0} = 1$, $\genfrac\lvert\rvert{0pt}{}{n} {0} = \genfrac\lvert\rvert{0pt}{}{0} {k} = 0$ and $\genfrac\lvert\rvert{0pt}{}{n} {k} = 0$ for $k > n$.
\end{defn}

\begin{thm}
For $n\geq k \geq 1$, the Ward-Lah numbers satisfy an explicit formula
\begin{align}
\label{wardlahexplicit}
\genfrac\lvert\rvert{0pt}{}{n} {k} = \frac{(n+k)!}{k!} \binom{n-1}{k-1}.
\end{align}
\end{thm}

\begin{proof}
From the definition \ref{wardLahdef},  (\ref{fallingfactiroalexplicitdef}) and (\ref{pnk1,1,1}), we get
\begin{align*}
\genfrac\lvert\rvert{0pt}{}{n} {k} &= (-1)^{k} (n+k)^{\underline{n}} P_{n}^{k} (1, 1, \cdots) \\
&= (-1)^{k} (n+k)^{\underline{n}} (-1)^{k} \binom{n-1}{k-1} \\
&= \frac{(n+k)!}{k!} \binom{n-1}{k-1}.
\end{align*}
\end{proof}

\begin{thm}
The Ward-Lah numbers satisfy
\begin{align}
\genfrac\lvert\rvert{0pt}{}{n} {k} &=  \sum_{m=0}^{k} (-1)^{m+k} \binom{n+k}{n+m} \binom{n+m-1}{m-1} \frac{(n+m)!}{m!} \nonumber \\
&=  \sum_{m=0}^{k} (-1)^{m+k} \binom{n+k}{n+m} \lah{n+m}{m} \label{wardlah-lah}.
\end{align}
\end{thm}

\begin{proof}
Simplifying both sides using explicit formulas (\ref{wardlahexplicit}), (\ref{lahexplicit}) and the explicit formula for binomial coefficients, we get 
\begin{align}
\label{proofsum}
\frac{(n+k)!(n-1)!}{k!(k-1)!(n-k)!} = \sum_{m=0}^{k} (-1)^{m+k} \frac{(n+k)!(n+m-1)!}{m!n!(m-1)!(k-m)!}.
\end{align}
Using Gosper's algorithm (\cite[pp. 75]{zeilbergerpetkovsekA=B}), we verify (\ref{proofsum}).  Therefore, Ward-Lah numbers coincide with the triangular array given in \cite[A357367]{SloaneOEIS}.
\end{proof}

\begin{rem}
Comparing (\ref{lahpartitiondef}) with (\ref{wardlahdefformula}) and because of (\ref{wardlah-lah}), we understand why we call these numbers Ward-Lah numbers. 
\end{rem}

\begin{cor}
Applying (\ref{wardlahexplicit}), we get some special values for Ward-Lah numbers.
\begin{align*}
\genfrac\lvert\rvert{0pt}{}{n} {1} = (n+1)! \\
\genfrac\lvert\rvert{0pt}{}{n} {n} = \frac{(2n)!}{n!} \\
\genfrac\lvert\rvert{0pt}{}{n} {n-1} = \frac{(2n-1)!}{(n-2)!}
\end{align*}
\end{cor}

\subsection{Triangular Recurrence Relations}

Here, we give triangular recurrence relations for Ward-Lah numbers. 

\begin{thm}
For $n \geq k \geq 1$, $k-1 \geq 1$, the Ward-Lah numbers satisfy the recurrence relation
\begin{align*}
\genfrac\lvert\rvert{0pt}{}{n} {k} = \frac{(n+k)(n-1)}{n} \left( \genfrac\lvert\rvert{0pt}{}{n-1} {k} + \frac{n+k-1}{k-1} \genfrac\lvert\rvert{0pt}{}{n-1} {k-1}\right).
\end{align*}
\end{thm}

\begin{proof}
Using the explicit formula (\ref{wardlahexplicit}), we get
\begin{align*}
\frac{(n+k)!(n-1)!}{k!(k-1)!(n-k)!} &= \frac{(n+k)(n-1)}{n} \cdot \frac{(n+k-1)!(n-2)!}{k!(k-1)!(n-k-1)!} \\ &+ \frac{(n+k)(n-1)}{n} \cdot \frac{(n+k-1)(n+k-2)!(n-2)!}{(k-1)(k-2)!(k-1)!(n-k)!}  \\
&= \frac{(n+k)(n-1)}{n} \cdot \frac{(n+k-1)!(n-2)!(k-1)!(n-k)!}{k!(n-k-1)!((k-1)!)^{2}(n-k)!} \\ &+  \frac{(n+k)(n-1)}{n} \cdot \frac{(n+k-1)!(n-2)!k!(n-k-1)!}{k!(n-k-1)!((k-1)!)^{2}(n-k)!}  \\
&=  \frac{(n+k)(n-1)}{n} \cdot \frac{kn!(n+k-1)!}{(n-1)(k!)^{2} (n-k)!} \\
&= \frac{(n+k)!(n-1)!}{k!(k-1)!(n-k)!}.
\end{align*}
\end{proof}

\begin{thm}
For $n \geq k \geq 1$, the following recurrence relation with integer coefficients holds
\begin{align*}
\genfrac\lvert\rvert{0pt}{}{n} {k} = 2(n+k-1) \genfrac\lvert\rvert{0pt}{}{n-1} {k-1} + (n+2k-1) \genfrac\lvert\rvert{0pt}{}{n-1} {k}.
\end{align*}
\end{thm}

\begin{proof}
Using (\ref{wardlahexplicit}), we get
\begin{align*}
\frac{(n+k)!(n-1)!}{k!(k-1)!(n-k)!} &= \frac{2(n+k-1)!(n-2)!}{(k-1)!(k-2)!(n-k)!} + \frac{(n+2k-1)(n+k-2)!(n-2)!}{k!(k-1)!(n-k-1)!} \\
&= \frac{k(n-1)!(n+k)!}{(k!)^{2} (n-k)!} \\
&= \frac{(n+k)!(n-1)!}{k!(k-1)!(n-k)!}. 
\end{align*}
\end{proof}

\subsection{Horizontal Recurrence Relation}
Next, we give a horizontal recurrence relation for Ward-Lah numbers.

\begin{thm}
For positive integers $n, k, m$ and $ n\geq k \geq 1$, $n-m \geq 1$, $n-m+k-j \geq 1$ the Ward-Lah numbers satisfy a horizontal recurrence relation
\begin{align}
\label{horizontalwardlah}
\genfrac\lvert\rvert{0pt}{}{n} {k} = \frac{(n+k)!}{k!} \sum_{j=0}^{m} \frac{(k-j)!}{(n-m+k-j)!} \binom{m}{j} \genfrac\lvert\rvert{0pt}{}{n-m} {k-j}.
\end{align}
\end{thm}

\begin{proof}
Using Vandermonde's identity and (\ref{wardlahexplicit}), we get
\begin{align*}
\genfrac\lvert\rvert{0pt}{}{n} {k} &=  \frac{(n+k)!}{k!} \binom{n-1}{k-1} \\
&= \frac{(n+k)!}{k!} \sum_{j=0}^{m} \binom{m}{j} \binom{n-m-1}{k-j-1}.
\end{align*}
Note that $\binom{n-m-1}{k-j-1} = \frac{(k-j)!}{(n-m+k-j)!} \genfrac\lvert\rvert{0pt}{}{n-m} {k-j}$. The result follows. 
\end{proof}

\begin{cor}
From horizontal recurrence relation for $m=1$ and $n \geq k \geq 1$, we get another triangular recurrence for Ward-Lah numbers
\begin{align*}
\genfrac\lvert\rvert{0pt}{}{n} {k} = (n+k) \left( \genfrac\lvert\rvert{0pt}{}{n-1} {k} + \frac{n+k-1}{k} \genfrac\lvert\rvert{0pt}{}{n-1} {k-1}\right).
\end{align*}
\end{cor}

\begin{proof}
Using (\ref{horizontalwardlah}) for $m=1$, we get
\begin{align*}
\genfrac\lvert\rvert{0pt}{}{n} {k} = \frac{k!(n+k)!}{k!(n+k-1)!} \genfrac\lvert\rvert{0pt}{}{n-1} {k} + \frac{(n+k)!(n-1)!}{(n+k-2)!k!} \genfrac\lvert\rvert{0pt}{}{n-1} {k-1}.
\end{align*}
The result follows. Note that we get the same result via explicit formula.
\end{proof}

\subsection{Recurrence Relation of Order 3}

Here, we present a recurrence relation of order 3 obtained by using Sister Celine's general algorithm \cite[pp. 59]{zeilbergerpetkovsekA=B}.

\begin{thm}
For $n \geq 2$ and $k \geq 1$, the Ward-Lah numbers satisfy
\begin{align*}
\genfrac\lvert\rvert{0pt}{}{n} {k} = 2(2n-1) \genfrac\lvert\rvert{0pt}{}{n-1} {k-1} - n(n-2) \genfrac\lvert\rvert{0pt}{}{n-2} {k} - (-2n+1)\genfrac\lvert\rvert{0pt}{}{n-1} {k}.
\end{align*}
\end{thm}

\begin{proof}
Using (\ref{wardlahexplicit}), we get the result. 
\end{proof}

\subsection{Exponential Generating Function}

We now give an exponential generating function for Ward-Lah numbers. 

\begin{thm}
For $n \geq k$, the Ward-Lah numbers satisfy an exponential generating function 
\begin{align*}
\sum_{n=k}^{\infty} \genfrac\lvert\rvert{0pt}{}{n-k} {k} \frac{x^{n}}{n!} = \frac{1}{k!} \left( \frac{x^{2k}}{(1-x)^{k}} \right).
\end{align*}
\end{thm}
\begin{proof}
Using (\ref{wardlahexplicit}) and the binomial series, we get
\begin{align*}
\sum_{n=k}^{\infty} \genfrac\lvert\rvert{0pt}{}{n-k} {k} \frac{x^{n}}{n!} &= \sum_{n=k}^{\infty} \frac{n!}{k!} \binom{n-k-1}{k-1} \frac{x^{n}}{n!} \\
&= \sum_{n=k}^{\infty} \frac{1}{k!} \binom{n-k-1}{k-1} x^{n} \\
&= \frac{x^{2k}}{k!} \sum_{n=k}^{\infty} \binom{n-k-1}{k-1} x^{n-2k} \\
&=  \frac{x^{2k}}{k!} \sum_{n=k}^{\infty} \binom{n-k-1}{n-2k} x^{n-2k} \\
&=  \frac{x^{2k}}{k!} \sum_{n=0}^{\infty} \binom{n-k-1+2k}{n-2k+2k} x^{n} \\
&=  \frac{x^{2k}}{k!} \sum_{n=0}^{\infty} \binom{n+k-1}{n} x^{n} \\
&=  \frac{x^{2k}}{k!} \sum_{n=0}^{\infty} \binom{-k}{n} (-x)^{n} \\
&=  \frac{x^{2k}}{k!} (1-x)^{-k} 
=   \frac{1}{k!} \left( \frac{x^{2k}}{(1-x)^{k}} \right).  \end{align*}
\end{proof}

\section{Varied Ward Numbers} \label{section3}

Here, we give recurrence relations for sequences defined by equations (\ref{t1defintro}) and (\ref{t2defintro}). We now use notations

\begin{align*}
T_{1}^{*}(n,k) = \left\uparrow
\begin{matrix}
n \\k
\end{matrix}
\right\downarrow^{*} \\
T_{2}^{*}(n,k) = \left\updownarrow
\begin{matrix}
n \\k
\end{matrix}
\right\updownarrow^{*}.
\end{align*} Comparing (\ref{t1defintro}) with (\ref{kPWard1def}) and (\ref{t2defintro}) with (\ref{kPWard2def}), we get

\begin{align}
\left\uparrow
\begin{matrix}
n \\k
\end{matrix}
\right\downarrow^{*} = \label{t1def}
\frac{(2n)!}{(n+k)^{\underline{n}}} \left\uparrow
\begin{matrix}
n \\k
\end{matrix}
\right\downarrow
\end{align}
and 
\begin{align}
\label{t2def} \left\updownarrow
\begin{matrix}
n \\k
\end{matrix}
\right\updownarrow^{*} = 
\frac{(2n)!}{(n+k)^{\underline{n}}} 	\left\updownarrow
\begin{matrix}
n \\k
\end{matrix}
\right\updownarrow.
\end{align}

\begin{rem}
Since $(n+k)^{\underline{n}}$ is the number of $k$-element variations of $n$ objects, we call these numbers \textit{varied Ward numbers of both kinds}. 
\end{rem}

\subsection{Triangular Recurrence Relation for Varied Ward Numbers of the First Kind}

\begin{thm}
For $n \geq k \geq 1$, the varied Ward numbers of the first kind satisfy the recurrence relation
\begin{align*}
\left\uparrow
\begin{matrix}
n \\k
\end{matrix}
\right\downarrow^{*} = \frac{2n(2n-1)}{n+k} \left( (n+k-1) \left\uparrow
\begin{matrix}
n-1 \\k
\end{matrix}
\right\downarrow^{*} + k \left\uparrow
\begin{matrix}
n-1 \\k-1
\end{matrix}
\right\downarrow^{*}  \right)
\end{align*}
with boundary conditions, given by (\ref{boundary1}), (\ref{boundary2}) and (\ref{boundary 2a}). 
\end{thm}

\begin{proof}
Using (\ref{t1def}), we write the recurrence relation for varied Ward numbers of the first kind in terms of Ward numbers of the first kind (\ref{ward1rec})
\begin{align*}
\left\uparrow
\begin{matrix}
n \\k
\end{matrix}
\right\downarrow^{*} = \frac{(2n)!k!}{(n+k)!} (n+k-1) \left\uparrow
\begin{matrix}
n-1 \\k
\end{matrix}
\right\downarrow +  \frac{(2n)!k!}{(n+k)!} (n+k-1) \left\uparrow
\begin{matrix}
n-1 \\k-1
\end{matrix}
\right\downarrow.
\end{align*}
Note that
\begin{align*}
\left\uparrow
\begin{matrix}
n-1 \\k
\end{matrix}
\right\downarrow^{*} = \frac{(2(n-1))!k!}{(n+k-1)!} \left\uparrow
\begin{matrix}
n-1 \\k
\end{matrix}
\right\downarrow
\end{align*}
and
\begin{align*}
\left\uparrow
\begin{matrix}
n-1 \\k-1
\end{matrix}
\right\downarrow^{*} = \frac{(2(n-1))!(k-1)!}{(n+k-2)!} \left\uparrow
\begin{matrix}
n-1 \\k-1
\end{matrix}
\right\downarrow.
\end{align*}
Thus, we get
\begin{align*}
\left\uparrow
\begin{matrix}
n \\k
\end{matrix}
\right\downarrow^{*} &= \frac{\frac{(2n)!k!}{(n+k)!}}{\frac{(2(n-1))!k!}{(n+k-1)!}} (n+k-1) \left\uparrow
\begin{matrix}
n-1 \\k
\end{matrix}
\right\downarrow^{*}  + \frac{\frac{(2n)!k!}{(n+k)!}}{\frac{(2(n-1))!(k-1)!}{(n+k-2)!}} (n+k-1) \left\uparrow
\begin{matrix}
n-1 \\k-1
\end{matrix}
\right\downarrow^{*} \\
&= \frac{2n(2n-1)}{n+k} (n+k-1) \left\uparrow
\begin{matrix}
n-1 \\k
\end{matrix}
\right\downarrow^{*}  +  \frac{2nk(2n-1)}{n+k} \left\uparrow
\begin{matrix}
n-1 \\k-1
\end{matrix}
\right\downarrow^{*} .
\end{align*}
By taking the common factor $\frac{2n(2n-1)}{n+k}$ out, we get the result. 
\end{proof}

\subsection{Triangular Recurrence Relation for Varied Ward Numbers of the Second Kind}

\begin{thm}
For $n \geq k \geq 1$, the varied Ward numbers of the second kind satisfy the recurrence relation
\begin{align*}
\left\updownarrow
\begin{matrix}
n \\k
\end{matrix}
\right\updownarrow^{*} = \frac{2nk(2n-1)}{n+k} \left( \left\updownarrow
\begin{matrix}
n-1 \\k
\end{matrix}
\right\updownarrow^{*} + \left\updownarrow
\begin{matrix}
n-1 \\k-1
\end{matrix}
\right\updownarrow^{*}  \right)
\end{align*}
with boundary conditions, given by (\ref{boundary1}), (\ref{boundary2}) and (\ref{boundary 2a}).
\end{thm}

\begin{proof}
Using (\ref{t2def}), we write the recurrence relation for varied Ward numbers of the second kind in terms of Ward numbers of the second kind (\ref{ward2rec})
\begin{align*}
\left\updownarrow
\begin{matrix}
n \\k
\end{matrix}
\right\updownarrow^{*} = \frac{(2n)!k!}{(n+k)!} k \left\updownarrow
\begin{matrix}
n-1 \\k
\end{matrix}
\right\updownarrow +  \frac{(2n)!k!}{(n+k)!} (n+k-1) \left\updownarrow
\begin{matrix}
n-1 \\k-1
\end{matrix}
\right\updownarrow.
\end{align*}
Note that
\begin{align*}
\left\updownarrow
\begin{matrix}
n-1 \\k
\end{matrix}
\right\updownarrow^{*} = \frac{(2(n-1))!k!}{(n+k-1)!} \left\updownarrow
\begin{matrix}
n-1 \\k
\end{matrix}
\right\updownarrow
\end{align*}
and
\begin{align*}
\left\updownarrow
\begin{matrix}
n-1 \\k-1
\end{matrix}
\right\updownarrow^{*} = \frac{(2(n-1))!(k-1)!}{(n+k-2)!} \left\updownarrow
\begin{matrix}
n-1 \\k-1
\end{matrix}
\right\updownarrow.
\end{align*}
Thus, we get
\begin{align*}
\left\updownarrow
\begin{matrix}
n \\k
\end{matrix}
\right\updownarrow^{*} &= \frac{\frac{(2n)!k!}{(n+k)!}}{\frac{(2(n-1))!k!}{(n+k-1)!}} k \left\updownarrow
\begin{matrix}
n-1 \\k
\end{matrix}
\right\updownarrow^{*}  + \frac{\frac{(2n)!k!}{(n+k)!}}{\frac{(2(n-1))!(k-1)!}{(n+k-2)!}} (n+k-1) \left\updownarrow
\begin{matrix}
n-1 \\k-1
\end{matrix}
\right\updownarrow^{*} \\
&= k \frac{2n(2n-1)}{n+k} \left\updownarrow
\begin{matrix}
n-1 \\k
\end{matrix}
\right\updownarrow^{*}  +  \frac{2nk(2n-1)}{n+k} \left\updownarrow
\begin{matrix}
n-1 \\k-1
\end{matrix}
\right\updownarrow^{*} .
\end{align*}
By taking the common factor $\frac{2nk(2n-1)}{n+k}$ out, we get the result. 
\end{proof}

\section{Varied Ward-Lah Numbers} \label{section4}

\subsection{Definition and Explicit Formula}

We now give recurrence relations and some identities for the sequence (\ref{t3def}), which we now call \textit{varied Ward-Lah numbers} $\genfrac\lvert\rvert{0pt}{}{n} {k}^{*} $. These numbers were introduced by Luschny in \cite{LuschnyGitHub}. This sequence is not yet entered in \cite{SloaneOEIS}. We now use notation $\genfrac\lvert\rvert{0pt}{}{n} {k}^{*} = T_{3}^{*}(n,k)$.

\begin{defn}
\label{variedwardLahdef}
The varied Ward-Lah numbers are for $n, k \in \mathbb{N}_{0}$ and $n\geq k$ defined via the Partition transformation as
\begin{align*}
\genfrac\lvert\rvert{0pt}{}{n} {k}^{*} = (-1)^{k} (2n)! P_{n}^{k} (1, 1, \cdots)
\end{align*}
with boundary conditions, given by (\ref{boundary1}), (\ref{boundary2}) and (\ref{boundary 2a}).
\end{defn}

\begin{thm}
For $n\geq k \geq 1$, the varied Ward-Lah numbers satisfy an explicit formula
\begin{align}
\label{variedwardlahexplicit}
\genfrac\lvert\rvert{0pt}{}{n} {k}^{*} = (2n)! \binom{n-1}{k-1}.
\end{align}
\end{thm}

\begin{proof}
From the definition \ref{variedwardLahdef}, and (\ref{pnk1,1,1}), we get
\begin{align*}
\genfrac\lvert\rvert{0pt}{}{n} {k}^{*} &= (-1)^{k} (2n)! P_{n}^{k} (1, 1, \cdots) \\
&= (-1)^{k} (2n)! (-1)^{k} \binom{n-1}{k-1} \\
&=  (2n)! \binom{n-1}{k-1}.
\end{align*}
\end{proof}

\begin{cor}
Applying (\ref{variedwardlahexplicit}), we get some special values for varied Ward-Lah numbers. 
\begin{align*}
\genfrac\lvert\rvert{0pt}{}{n} {1}^{*} = \genfrac\lvert\rvert{0pt}{}{n} {n}^{*} = (2n)! \\
\genfrac\lvert\rvert{0pt}{}{n} {n-1}^{*} = (n-1)(2n)!
\end{align*}
\end{cor}

\begin{rem}
Since varied Ward-Lah numbers can be written in terms of Ward-Lah numbers as 
\begin{align*}
\genfrac\lvert\rvert{0pt}{}{n} {k}^{*} = \frac{(2n)!}{(n+k)^{\underline{n}}} \genfrac\lvert\rvert{0pt}{}{n} {k}
\end{align*} and since $(n+k)^{\underline{n}}$ is the number of $k$-element variations of $n$ objects we understand why we call these numbers varied Ward-Lah numbers.
\end{rem}

\subsection{Triangular Recurrence Relation}

Here, we give a triangular recurrence relation for varied Ward-Lah numbers. 

\begin{thm}
For $n \geq k \geq 1$, the varied Ward-Lah numbers satisfy the recurrence relation
\begin{align}
\label{variedWardLahtriangular}
\genfrac\lvert\rvert{0pt}{}{n} {k}^{*} = 2n(2n-1) \left( \genfrac\lvert\rvert{0pt}{}{n-1} {k}^{*} + \genfrac\lvert\rvert{0pt}{}{n-1} {k-1}^{*} \right).
\end{align}
\end{thm}

\begin{proof}
Using the explicit formula (\ref{variedwardlahexplicit}), we get
\begin{align*}
\frac{(2n)!(n-1)!}{(k-1)!(n-k)!} &= 2n(2n-1) \left( \frac{(2(n-1))!(n-2)!}{(k-1)!(n-k-1)!} + \frac{(2(n-1))!(n-2)!}{(k-2)!(n-k)!} \right)  \\
&= 2n(2n-1)  \frac{(2(n-1))!(n-2)!(k-2)!(n-k)!}{(k-1)!(k-2)!(n-k-1)!(n-k)!}  \\ &+  2n(2n-1)  \frac{(2(n-1))!(n-2)!(k-1)!(n-k-1)!}{(k-1)!(k-2)!(n-k-1)!(n-k)!}  \\
&=2n(2n-1)  \frac{((n-k)+(k-1))((2(n-1))!(n-2)!(k-2)!(n-k-1))}{(k-1)!(k-2)!(n-k-1)!(n-k)!}  \\
&= 2n(2n-1) \frac{(2(n-1))!(n-1)!}{(k-1)!(n-k)!} \\
&= \frac{(2n)!(n-1)!}{(k-1)!(n-k)!}.
\end{align*}
\end{proof}

\begin{rem}
By using Sister Celine's general algorithm \cite[pp. 59]{zeilbergerpetkovsekA=B}, we obtain again the triangular recurrence (\ref{variedWardLahtriangular}). 
\end{rem}

\subsection{Horizontal Recurrence Relation}
Next, we give a horizontal recurrence relation for varied Ward-Lah numbers.

\begin{thm}
For positive integers $n, k, m$ and $ n\geq k \geq 1$, $n-m \geq 1$, the varied Ward-Lah numbers satisfy a horizontal recurrence relation
\begin{align*}
\genfrac\lvert\rvert{0pt}{}{n} {k}^{*} = (2n)! \sum_{j=0}^{m} \frac{1}{(2(n-m))!} \binom{m}{j} \genfrac\lvert\rvert{0pt}{}{n-m} {k-j}^{*}.
\end{align*}
\end{thm}

\begin{proof}
Using Vandermonde's identity and (\ref{variedwardlahexplicit}), we get
\begin{align*}
\genfrac\lvert\rvert{0pt}{}{n} {k}^{*} &=  (2n)! \binom{n-1}{k-1} \\
&= (2n)! \sum_{j=0}^{m} \binom{m}{j} \binom{n-m-1}{k-j-1}.
\end{align*}
Note that $\binom{n-m-1}{k-j-1} = \frac{1}{(2(n-m))!} \genfrac\lvert\rvert{0pt}{}{n-m} {k-j}^{*}$. The result follows. 
\end{proof}

\begin{rem}
From horizontal recurrence relation for $m=1$ we get (\ref{variedWardLahtriangular}).
\end{rem}

\subsection{Generating Function}

Now, we give a generating function for varied Ward-Lah numbers. 

\begin{thm}
The generating function for the varied Ward-Lah numbers is
\begin{align*}
\sum_{n=k}^{\infty} \genfrac\lvert\rvert{0pt}{}{n} {k}^{*} \frac{x^{n}}{(2n)!} =  \left( \frac{x}{1-x} \right)^{k}.
\end{align*}
\end{thm}
\begin{proof}
Using (\ref{variedwardlahexplicit}) and the binomial series, we get
\begin{align*} 
\sum_{n=k}^{\infty} \genfrac\lvert\rvert{0pt}{}{n} {k}^{*} \frac{x^{n}}{(2n)!} &= \sum_{n=k}^{\infty} (2n)! \binom{n-1}{k-1} \frac{x^{n}}{(2n)!} \\
&= \sum_{n=k}^{\infty} \binom{n-1}{k-1} x^{n} \\
&= x^{k} \sum_{n=k}^{\infty} \binom{n-1}{k-1} x^{n-k} \\
&= x^{k} \sum_{n=0}^{\infty} \binom{n+k-1}{n} x^{n} \\
&= x^{k} \sum_{n=0}^{\infty} \binom{-k}{n} (-x)^{n} \\
&= x^{k} (1-x)^{-k} = \frac{x^{k}}{(1-x)^{k}} = \left( \frac{x}{1-x} \right)^{k}.
\end{align*}
\end{proof}

\subsection{Relation With Lah Numbers}

Luschny gave identities connecting varied Ward number of both kinds and Stirling numbers of both kinds (see \cite[Formula 14, Formula 15]{LuschnyGitHub}). For more information about Stirling numbers, see \cite{petkovsekpisanski}.

Motivated by Luschny's identities, we give an analogous identity for Lah and varied Ward-Lah numbers. 
\begin{thm}
For $n \geq k$, the following identity holds
\begin{align*}
(n-k+1)^{\overline{n-k}} \lah{n}{k} = \binom{n}{k} \sum_{j=0}^{k} \binom{k}{j} \genfrac\lvert\rvert{0pt}{}{n-k} {j}^{*},
\end{align*}   
where $(n-k+1)^{\overline{n-k}}$ is the rising factorial, defined by $x^{\overline{n}} = \frac{(x+n-1)!}{(x-1)!}$.
\end{thm}  

\begin{proof}
Simplifying the left-hand side using explicit formula (\ref{lahexplicit}), explicit formula for binomial coefficients and explicit formula for rising factorials, we get
\begin{align*}
(n-k+1)^{\overline{n-k}} \lah{n}{k} &= \frac{n!(2n-2k)!(n-1)!}{k!(k-1)!((n-k)!)^{2}} \\
&= (2(n-k))! \binom{n}{k} \binom{n-1}{k-1}.
\end{align*}
Using Vandermonde's identity, we get
\begin{align*}
(2(n-k))! \binom{n}{k} \binom{n-1}{k-1} =  (2(n-k))! \binom{n}{k} \sum_{j=0}^{m} \binom{m}{j} \binom{n-m-1}{k-j-1}.
\end{align*}
Setting $m=k$ gives
\begin{align*}
(2(n-k))! \binom{n}{k} \binom{n-1}{k-1} &= \binom{n}{k} \sum_{j=0}^{k} \binom{k}{j} (2(n-k))! \binom{n-k-1}{k-j-1} \\
&= \binom{n}{k} \sum_{j=0}^{k} \binom{k}{j} \genfrac\lvert\rvert{0pt}{}{n-k} {k-j}^{*}.
\end{align*}
Note that 
\begin{align*}
\sum_{j=0}^{k}  \genfrac\lvert\rvert{0pt}{}{n-k} {k-j}^{*} = \sum_{j=0}^{k}  \genfrac\lvert\rvert{0pt}{}{n-k} {j}^{*}.
\end{align*}
Thus, we get the result. 
\end{proof}

\section{Binomial Ward Numbers} \label{section5}

In this section, we give recurrence relations for sequences defined by equations (\ref{t1defcircintro}) and (\ref{t2defcircintro}). We also give two conjectured relations between binomial Ward numbers and central Stirling numbers. We now use notations

\begin{align*}
T_{1}^{\circ}(n,k) = \left\uparrow
\begin{matrix}
n \\k
\end{matrix}
\right\downarrow^{\circ} \\
T_{2}^{\circ}(n,k) = \left\updownarrow
\begin{matrix}
n \\k
\end{matrix}
\right\updownarrow^{\circ}.
\end{align*} Comparing (\ref{t1defcircintro}) with (\ref{kPWard1def}) and (\ref{t2defcircintro}) with (\ref{kPWard2def}), we get

\begin{align}
\label{t1defcirc}
\left\uparrow
\begin{matrix}
n \\k
\end{matrix}
\right\downarrow^{\circ} = \binom{2n}{n+k} \left\uparrow
\begin{matrix}
n \\k
\end{matrix}
\right\downarrow
\end{align}

\begin{align}
\label{t2defcirc}
\left\updownarrow
\begin{matrix}
n \\k
\end{matrix}
\right\updownarrow^{\circ} = \binom{2n}{n+k} \left\updownarrow
\begin{matrix}
n \\k
\end{matrix}
\right\updownarrow.
\end{align}

\begin{rem}
Since $\binom{2n}{n+k}$ is the binomial coefficient, we call these numbers binomial Ward numbers. 
\end{rem}

\subsection{Triangular Recurrence Relations}

We now study triangular recurrence relations for binomial Ward numbers of both kinds. 

\begin{thm}
For $ n \geq k \geq 1$ and $n-k \geq 1$, the binomial Ward numbers of the first kind satisfy the recurrence relation
\begin{align*}
\left\uparrow
\begin{matrix}
n \\k
\end{matrix}
\right\downarrow^{\circ} = \frac{2n(2n-1)}{n+k} \left( \frac{n+k-1}{n-k} \left\uparrow
\begin{matrix}
n-1 \\k
\end{matrix}
\right\downarrow^{\circ} +  \left\uparrow
\begin{matrix}
n-1 \\k-1
\end{matrix}
\right\downarrow^{\circ}  \right)
\end{align*}
with boundary conditions, given by (\ref{boundary3}), (\ref{boundary4}) and (\ref{boundary 4a}).
\end{thm}

\begin{proof}
Using (\ref{t1defcirc}), we write the recurrence relation for binomial Ward numbers of the first kind in terms of Ward numbers of the first kind (\ref{ward1rec})
\begin{align*}
\left\uparrow
\begin{matrix}
n \\k
\end{matrix}
\right\downarrow^{\circ} = \frac{(2n)!}{(n+k)!(n-k)!} (n+k-1) \left\uparrow
\begin{matrix}
n-1 \\k
\end{matrix}
\right\downarrow +  \frac{(2n)!}{(n+k)!(n-k)!} (n+k-1) \left\uparrow
\begin{matrix}
n-1 \\k-1
\end{matrix}
\right\downarrow.
\end{align*}
Note that
\begin{align*}
\left\uparrow
\begin{matrix}
n-1 \\k
\end{matrix}
\right\downarrow^{\circ} = \frac{(2(n-1))!}{(n+k-1)!(n-k-1)!} \left\uparrow
\begin{matrix}
n-1 \\k
\end{matrix}
\right\downarrow
\end{align*}
and
\begin{align*}
\left\uparrow
\begin{matrix}
n-1 \\k-1
\end{matrix}
\right\downarrow^{\circ} = \frac{(2(n-1))!}{(n+k-2)!(n-k)!} \left\uparrow
\begin{matrix}
n-1 \\k-1
\end{matrix}
\right\downarrow.
\end{align*}
Thus, we get
\begin{align*}
\left\uparrow
\begin{matrix}
n \\k
\end{matrix}
\right\downarrow^{\circ} &= \frac{\frac{(2n)!}{(n+k)!(n-k)!}}{\frac{(2(n-1))!}{(n+k-1)!(n-k-1)!}} (n+k-1) \left\uparrow
\begin{matrix}
n-1 \\k
\end{matrix}
\right\downarrow^{\circ}  + \frac{\frac{(2n)!}{(n+k)!(n-k)!}}{\frac{(2(n-1))!}{(n+k-2)!(n-k)!}} (n+k-1) \left\uparrow
\begin{matrix}
n-1 \\k-1
\end{matrix}
\right\downarrow^{\circ} \\
&= \frac{2n(2n-1)}{(n+k)(n-k)} (n+k-1) \left\uparrow
\begin{matrix}
n-1 \\k
\end{matrix}
\right\downarrow^{\circ}  +  \frac{2n(2n-1)}{(n+k)} \left\uparrow
\begin{matrix}
n-1 \\k-1
\end{matrix}
\right\downarrow^{\circ} .
\end{align*}
By taking the common factor $\frac{2n(2n-1)}{n+k}$ out, we get the result. 
\end{proof}

\begin{thm}
For $ n \geq k \geq 1$ and $n-k \geq 1$, the binomial Ward numbers of the second kind satisfy the recurrence relation
\begin{align*}
\left\updownarrow
\begin{matrix}
n \\k
\end{matrix}
\right\updownarrow^{\circ} = \frac{2n(2n-1)}{n+k} \left( \frac{k}{n-k}  \left\updownarrow
\begin{matrix}
n-1 \\k
\end{matrix}
\right\updownarrow^{\circ} + \left\updownarrow
\begin{matrix}
n-1 \\k-1
\end{matrix}
\right\updownarrow^{\circ}  \right)
\end{align*}
with boundary conditions, given by (\ref{boundary3}), (\ref{boundary4}) and (\ref{boundary 4a}).
\end{thm}

\begin{proof}
Using (\ref{t2defcirc}), we write the recurrence relation for binomial Ward numbers of the second kind in terms of Ward numbers of the second kind (\ref{ward2rec})
\begin{align*}
\left\updownarrow
\begin{matrix}
n \\k
\end{matrix}
\right\updownarrow^{\circ} = \frac{(2n)!}{(n+k)!(n-k)!} k \left\updownarrow
\begin{matrix}
n-1 \\k
\end{matrix}
\right\updownarrow +  \frac{(2n)!}{(n+k)!(n-k)!} (n+k-1) \left\updownarrow
\begin{matrix}
n-1 \\k-1
\end{matrix}
\right\updownarrow.
\end{align*}
Note that
\begin{align*}
\left\updownarrow
\begin{matrix}
n-1 \\k
\end{matrix}
\right\updownarrow^{\circ} = \frac{(2(n-1))!}{(n+k-1)!(n-k-1)!} \left\updownarrow
\begin{matrix}
n-1 \\k
\end{matrix}
\right\updownarrow
\end{align*}
and
\begin{align*}
\left\updownarrow
\begin{matrix}
n-1 \\k-1
\end{matrix}
\right\updownarrow^{\circ} = \frac{(2(n-1))!}{(n+k-2)!(n-k)!} \left\updownarrow
\begin{matrix}
n-1 \\k-1
\end{matrix}
\right\updownarrow.
\end{align*}
Thus, we get
\begin{align*}
\left\updownarrow
\begin{matrix}
n \\k
\end{matrix}
\right\updownarrow^{\circ} &= \frac{\frac{(2n)!}{(n+k)!(n-k)!}}{\frac{(2(n-1))!}{(n+k-1)!(n-k-1)!}} k \left\updownarrow
\begin{matrix}
n-1 \\k
\end{matrix}
\right\updownarrow^{\circ}  + \frac{\frac{(2n)!}{(n+k)!(n-k)!}}{\frac{(2(n-1))!}{(n+k-2)!(n-k)!}} (n+k-1) \left\updownarrow
\begin{matrix}
n-1 \\k-1
\end{matrix}
\right\updownarrow^{\circ} \\
&= k \frac{2n(2n-1)}{(n+k)(n-k)} \left\updownarrow
\begin{matrix}
n-1 \\k
\end{matrix}
\right\updownarrow^{\circ}  +  \frac{2n(2n-1)}{(n+k)} \left\updownarrow
\begin{matrix}
n-1 \\k-1
\end{matrix}
\right\updownarrow^{\circ} .
\end{align*}
By taking the common factor $\frac{2n(2n-1)}{n+k}$ out, we get the result. 
\end{proof}

\subsection{Conjectured Relations With Central Stirling Numbers}

Here, we give two conjectured relations between binomial Ward numbers and central Stirling numbers based on experimental evidence. The central Stirling numbers of both kinds are entered in \cite[A007820, A187646]{SloaneOEIS}. $\stonet{2n}{n}$ denotes the central Stirling numbers of the first kind and $\sttwot{2n}{n}$ denotes the central Stirling numbers of the second kind, where $\stonet{n}{k}$ and $\sttwot{n}{k}$ denote unsigned Stirling numbers of both kinds (cf. \cite{petkovsekpisanski}). 

\begin{conj}
For $n \geq k$, the binomial Ward numbers satisfy 
\begin{align*}
\sum_{k=0}^{n} \left\uparrow
\begin{matrix}
n \\k
\end{matrix}
\right\downarrow^{\circ} = \stone{2n}{n} \\
\sum_{k=0}^{n} \left\updownarrow
\begin{matrix}
n \\k
\end{matrix}
\right\updownarrow^{\circ} = \sttwo{2n}{n}.
\end{align*}
\end{conj}

\section{Binomial Ward-Lah Numbers} \label{section6}

Finally, we analogously define the \textit{binomial Ward-Lah numbers}.

\begin{defn}
\label{binomialwardLahdef}
The binomial Ward-Lah numbers are for $n, k \in \mathbb{N}_{0}$ and $n-k \geq 1$ defined via the Partition transformation as
\begin{align*}
\genfrac\lvert\rvert{0pt}{}{n} {k}^{\circ} = (-1)^{k} \frac{(2n)!}{k!(n-k)!} P_{n}^{k} (1, 1, \cdots)
\end{align*}
with
\begin{align*}
\genfrac\lvert\rvert{0pt}{}{n} {0}^{\circ} = \genfrac\lvert\rvert{0pt}{}{0} {k}^{\circ} = 0
\end{align*}
and
\begin{align*}
\genfrac\lvert\rvert{0pt}{}{0} {0}^{\circ} = 1.
\end{align*}
Also, $\genfrac\lvert\rvert{0pt}{}{n} {k}^{\circ} = 0$ for $k > n$. 
\end{defn}

\begin{thm}
For $n\geq k \geq 1$ and $n-k \geq 1$, the binomial Ward-Lah numbers satisfy an explicit formula
\begin{align}
\label{binomialwardlahexplicit}
\genfrac\lvert\rvert{0pt}{}{n} {k}^{\circ} = \frac{(2n)!}{k!(n-k)!} \binom{n-1}{k-1}.
\end{align}
\end{thm}

\begin{proof}
From the definition \ref{binomialwardLahdef} and (\ref{pnk1,1,1}), we get
\begin{align*}
\genfrac\lvert\rvert{0pt}{}{n} {k}^{\circ} &= (-1)^{k} \frac{(2n)!}{k!(n-k)!} P_{n}^{k} (1, 1, \cdots) \\
&= (-1)^{k} \frac{(2n)!}{k!(n-k)!} (-1)^{k} \binom{n-1}{k-1} \\
&=  \frac{(2n)!}{k!(n-k)!} \binom{n-1}{k-1}.
\end{align*}
\end{proof}

\begin{cor}
Applying (\ref{binomialwardlahexplicit}), we get some special values for binomial
Ward-Lah numbers. 
\begin{align*}
\genfrac\lvert\rvert{0pt}{}{n} {1}^{\circ} = \frac{(2n)!}{(n-1)!} \\ \genfrac\lvert\rvert{0pt}{}{n} {n}^{\circ} = \frac{(2n)!}{n!}  \\
\genfrac\lvert\rvert{0pt}{}{n} {n-1}^{\circ} = \frac{(2n)!}{(n-2)!}
\end{align*}
\end{cor}

\begin{rem}
Since binomial Ward-Lah numbers can be written in terms of Ward-Lah numbers as 
\begin{align*}
\genfrac\lvert\rvert{0pt}{}{n} {k}^{\circ} = \binom{2n}{n+k} \genfrac\lvert\rvert{0pt}{}{n} {k}
\end{align*} and since $\binom{2n}{n+k}$ is the binomial coefficient we understand why we call these numbers binomial Ward-Lah numbers.
\end{rem}

\subsection{Triangular Recurrence Relation}

Here, we give a triangular recurrence relation for binomial Ward-Lah numbers. 

\begin{thm}
For $n \geq k \geq 1$, $n-k \geq 1$, the binomial Ward-Lah numbers satisfy the recurrence relation
\begin{align}
\label{binomialWardLahtriangular}
\genfrac\lvert\rvert{0pt}{}{n} {k}^{\circ} = 2n(2n-1) \left( \frac{1}{n-k} \genfrac\lvert\rvert{0pt}{}{n-1} {k}^{\circ} + \frac{1}{k} \genfrac\lvert\rvert{0pt}{}{n-1} {k-1}^{\circ} \right) .
\end{align}
\end{thm}

\begin{proof}
Using the explicit formula (\ref{binomialwardlahexplicit}), we get
\begin{align*}
\frac{(2n)!(n-1)!}{k!(k-1)!((n-k)!)^{2}} &=  2n(2n-1) \cdot \frac{(2(n-1))! (n-2)!}{(n-k)k!(k-1)!((n-k-1)!)^{2}}  \\ &+ 2n(2n-1) \cdot \frac{(2(n-1))!(n-2)!}{k(k-1)!(k-2)!((n-k)!)^{2}}  \\
&= 2n(2n-1) \cdot \frac{(2(n-1))!(n-2)!(k-2)!(n-k)!}{k!(k-1)!(k-2)!((n-k)!)^{2} (n-k-1)!}  \\ &+ 2n(2n-1)   \cdot \frac{(2(n-1))!(n-2)!(k-1)!(n-k-1)!}{k!(k-1)!(k-2)!((n-k)!)^{2} (n-k-1)!}  \\
&= 2n(2n-1) \frac{k(n-1)!(2(n-1))!}{(k!)^{2} ((n-k)!)^{2}} \\
&= 	\frac{(2n)!(n-1)!}{k!(k-1)!((n-k)!)^{2}}.
\end{align*}
\end{proof}

\subsection{Horizontal Recurrence Relation}
Now, we give a horizontal recurrence relation for binomial Ward-Lah numbers.

\begin{thm}
For positive integers $n, k, m$ and $n\geq k \geq 1$, $n-m \geq 1$, $n-k \geq 1$, the binomial Ward-Lah numbers satisfy a horizontal recurrence relation
\begin{align*}
\genfrac\lvert\rvert{0pt}{}{n} {k}^{\circ} = \frac{(2n)!}{k!(n-k)!} \sum_{j=0}^{m} \frac{(k-j)!(n-m-k+j)!}{(2(n-m))!} \binom{m}{j} \genfrac\lvert\rvert{0pt}{}{n-m} {k-j}^{\circ}.
\end{align*}
\end{thm}

\begin{proof}
Using Vandermonde's identity and (\ref{binomialwardlahexplicit}), we get
\begin{align*}
\genfrac\lvert\rvert{0pt}{}{n} {k}^{\circ} &=  \frac{(2n)!}{k!(n-k)!} \binom{n-1}{k-1} \\
&= \frac{(2n)!}{k!(n-k)!}  \sum_{j=0}^{m} \binom{m}{j} \binom{n-m-1}{k-j-1}.
\end{align*}
Note that $\binom{n-m-1}{k-j-1} = \frac{(k-j)!(n-m-k+j)!}{(2(n-m))!} \genfrac\lvert\rvert{0pt}{}{n-m} {k-j}^{\circ}$. The result follows. 
\end{proof}

\begin{rem}
From horizontal recurrence relation for $m=1$ we get (\ref{binomialWardLahtriangular}).
\end{rem}

\subsection{Recurrence Relation of Order 5}

Using sister Celine's general algorithm, we get another recurrence relation for binomial Ward numbers of order 5.

\begin{thm}
For $n, k \geq 2$ and $2n-3 \geq 0$, the following recurrence relation holds
\begin{align*}
\genfrac\lvert\rvert{0pt}{}{n} {k}^{\circ} &= \frac{-4(n-2)(2n-1)^{2}}{n} \left(\genfrac\lvert\rvert{0pt}{}{n-2} {k-2}^{\circ} -2 \genfrac\lvert\rvert{0pt}{}{n-2} {k-1}^{\circ} + \genfrac\lvert\rvert{0pt}{}{n-2} {k}^{\circ} \right) \\ &+ \frac{4(2n-1)}{n(2n-3)} \left((2(n-1)^{2}-1) \genfrac\lvert\rvert{0pt}{}{n-1} {k-1}^{\circ} + 2(n-1)^{2} \genfrac\lvert\rvert{0pt}{}{n-1} {k}^{\circ} \right). 
\end{align*}
\end{thm}

\begin{proof}
Using explicit formula (\ref{binomialwardlahexplicit}), we get the result. 
\end{proof}

\subsection{Relation With Central Lah Numbers}

Here, we give a relation connecting binomial Ward-Lah numbers and central Lah numbers $\laht{2n}{n}$. Central Lah numbers are entered in \cite[A187535]{SloaneOEIS}.

\begin{thm}
For $n-k \geq 1$, the binomial Ward-Lah numbers satisfy
\begin{align*}
\sum_{k=0}^{n} \genfrac\lvert\rvert{0pt}{}{n} {k}^{\circ} = \lah{2n}{n}.
\end{align*}
\end{thm}

\begin{proof}
Using (\ref{binomialwardlahexplicit}) and Gosper's algorithm (\cite[pp. 75]{zeilbergerpetkovsekA=B}), we get
\begin{align*}
\sum_{k=0}^{n} \frac{(2n)!}{k!(n-k)!} \binom{n-1}{k-1} = \frac{(2n)!}{n!} \binom{2n-1}{n-1}.
\end{align*}
Note that $\frac{(2n)!}{n!} \binom{2n-1}{n-1} = \laht{2n}{n}$. The result follows. 
\end{proof}

\section*{Acknowledgement}

The author is grateful to his parents for their support.

\end{document}